\documentclass[reqno,english]{smfart}
\usepackage[english]{babel}
\usepackage{amscd,smfthm}
\NumberTheoremsIn{subsection}

\author{Yves Laurent}
\address{ Institut Fourier Math\'ematiques \\
          UMR 5582  CNRS/UJF\\
                 BP 74\\
         38402 St Martin d'H\`eres Cedex}
\email{Yves.Laurent@ujf-grenoble.fr} \urladdr{http://www-fourier.ujf-grenoble.fr/~laurenty}
\title{\protect\boldmath$b$\protect\unboldmath-functions and integrable solutions of holonomic D-module}
\alttitle{\protect\boldmath$b$\protect\unboldmath-fonctions et solutions int{\'e}grables des modules
holonomes}

 \DeclareMathOperator{\ad}{ad}

\begin{document}
\frontmatter
\begin{abstract}
A famous theorem of Harish-Chandra shows that all invariant eigendistributions on a semisimple Lie
group are locally integrable functions. We give here an algebraic version of this theorem in terms
of polynomials associated with a holonomic $\mathcal D$-module.
\end{abstract}
\begin{altabstract}
Un c\'el\`ebre th\'eor\`eme de Harish-Chandra montre que les distributions invariantes propres sur
un groupe de Lie semisimple sont des fonctions localement int\'egrables. Nous donnons ici une
version alg\'ebrique de ce th\'eor\`eme en termes de polyn\^omes associ\'es {\`a} un $\mathcal D$-module
holonome.
\end{altabstract}

\keywords{Vanishing cycles, D-modules, symmetric pair, b-function}
\subjclass{35A27,35D10, 17B15}
\dedicatory{\`A Jean-Pierre Ramis, \`a l'occasion de son 60e anniversaire.}
\maketitle

\bibliographystyle{smfplain}

\newcommand{\x}{\times}
\newcommand{\abs}[1]{\lvert#1\rvert}
\newcommand{\ensemble}[2]{\{\,#1\mid#2\,\}}

\newcommand{\C}{{\mathbb C}}
\renewcommand{\L}{{\mathbb L}}
\newcommand{\N}{{\mathbb N}}
\newcommand{\R}{{\mathbb R}}
\newcommand{\Z}{{\mathbb Z}}

\newcommand{\CA}{\mathcal A}
\newcommand{\CaD}{\mathcal D}
\newcommand{\CF}{\mathcal F}
\newcommand{\CI}{\mathcal I}
\newcommand{\CK}{\mathcal K}
\newcommand{\CM}{\mathcal M}
\newcommand{\CN}{\mathcal N}
\newcommand{\CO}{\mathcal O}

\newcommand{\kN}{\mathfrak N}
\newcommand{\kO}{\mathfrak O}
\newcommand{\kc}{\mathfrak c}
\newcommand{\kg}{\mathfrak g}
\newcommand{\kh}{\mathfrak h}
\newcommand{\kl}{\mathfrak l}
\newcommand{\kq}{\mathfrak q}
\newcommand{\ks}{\mathfrak s}
\newcommand{\kw}{\mathfrak W}

\newcommand{\kgp}{\kg_P}
\newcommand{\kqp}{\kq_P}
\newcommand{\khp}{\kh_P}
\newcommand{\khpo}{\khp^\bot}
\newcommand{\khpp}{(\khpo)'}
\newcommand{\sld}{\ks\kl_2}

\newcommand{\vespa}{\vspace{1em}}

\newcommand{\TX}{{T^*X}}
\newcommand{\OX}{{{\mathcal O_X}}}
\newcommand{\DX}{{\mathcal D_X}}
\newcommand{\DcXY}{\widehat{\mathcal D}_{X|Y}}
\newcommand{\DcXYp}{\widehat{\mathcal D}_{X'|Y'}}
\newcommand{\DXp}{{\mathcal D_{X'}}}
\newcommand{\DXpX}{{\mathcal D_{X'\rightarrow X}}}
\newcommand{\DGV}{{\mathcal D_{\kg\rightarrow V}}}
\newcommand{\DXXp}{{\mathcal D_{X\leftarrow X'}}}

\newcommand{\scal}[2]{\langle #1,#2\rangle}

\def \projlim#1{\underset{#1}{\varprojlim}}
\def  \injlim#1{\underset{#1}{\varinjlim}}

\newcommand{\ga}{\alpha}
\newcommand{\gb}{\beta}
\newcommand{\gd}{\delta}           \newcommand{\gD}{\Delta}
\renewcommand{\gg}{\gamma}
\newcommand{\gh}{\eta}
\newcommand{\gk}{\kappa}
\newcommand{\gth}{\theta}          \newcommand{\gTh}{\Theta}
\newcommand{\gvt}{\vartheta}
\newcommand{\gl}{\lambda}          \newcommand{\gL}{\Lambda}
\newcommand{\gn}{\nu}
\newcommand{\gx}{\xi}              \newcommand{\gX}{\Xi}
\newcommand{\gp}{\pi}              \newcommand{\gP}{\Pi}
\newcommand{\gro}{\varrho}    
\newcommand{\gs}{\sigma}           \newcommand{\gS}{\Sigma}
\newcommand{\gt}{\tau}
\newcommand{\gf}{\varphi}          \newcommand{\gF}{\Phi}
\newcommand{\gy}{\psi}             \newcommand{\gY}{\Psi}
\newcommand{\go}{\omega}           \newcommand{\gO}{\Omega}

\newcommand{\exacte}[3]{0\xrightarrow{\ \ }{#1}\xrightarrow{\ \ }{#2}
 \xrightarrow{\ \ }{#3} \xrightarrow{\ \ }0}

\newcommand{\Dg}{\CaD_\kg}
\newcommand{\Dgg}{\CaD_\kg^G}
\newcommand{\Dh}{\CaD_\kh}
\newcommand{\DV}{\CaD_W}
\newcommand{\MF}{\CM_F}
\newcommand{\CHM}{{Ch(\CM)}}
\newcommand{\MFl}{\CM_\gl^F}
\newcommand{\IF}{\CI_F}
\newcommand{\Nl}{\CN_F}
\newcommand{\Mlh}{\CM_F^\kh}
\newcommand{\Nlh}{\CN_F^\kh}
\newcommand{\NO}{\CM_F^0}
\newcommand{\kgrs}{\kg_{rs}}
\newcommand{\gvtp}{{\gvt_P}}


\frontmatter
\section*{Introduction}

Let $G_\R$ be a real semisimple Lie group and $\kg_\R$ be its Lie algebra. An
\textsl{invariant eigendistribution} $T$ on $G_\R$ is a distribution which is invariant
under conjugation by elements of $G_\R$ and is an eigenvector of every bi-invariant
differential operator on $G_\R$. The main examples of such distributions are the
characters of irreducible representations of $G_\R$. A famous theorem of Harish-Chandra
sets that all invariant eigendistributions are $L^1_{loc}$-functions on $G_\R$
\cite{HC1}. After transfer to the Lie algebra by the exponential map, such a distribution
satisfies a system of partial differential equation.

In the language of $\CaD$-modules, these equations define a holonomic $\CaD$-modules on the
complexified Lie algebra $\kg$. We call this module the Hotta-Kashiwara module as it has been
defined and studied first in \cite{HOTTA}. In \cite{SEKI}, J. Sekiguchi extended these results to
symmetric pairs. He proved in particular that a condition on the symmetric pair is needed to extend
Harish-Chandra theorem. In several papers, Levasseur and Stafford \cite{LEVASS} \cite{LEVASS3}
\cite{LEVASS2} gave an algebraic proof of the main part of Harish-Chandra theorem.

In \cite{GL}, we defined a class of holonomic $\CaD$-modules, which we called \emph{tame
$\CaD$-modules}. These $\CaD$-modules have no quotients supported by a hypersurface and
their distribution solution are locally integrable. We proved in particular that the
Hotta-Kashiwara module is tame, recovering Harish-Chandra theorem. The definition of tame
is a condition on the roots of the $b$-functions which are polynomials attached to the
$\CaD$-module and a stratification of the base space. However, the proof of the fact that
the Hotta-Kashiwara module is tame involved some non algebraic vector fields.

The first aim of this paper is to give a completely algebraic version of Harish-Chandra
theorem. We give a slightly different definition of tame and an algebraic proof of the
fact that the Hotta-Kashiwara module is tame. This proof is different from the proof of
\cite{GL} and gives more precise results on the roots of the $b$-functions. However our
first proof was still valid in the case of symmetric pairs while the present proof uses a
morphism of Harish-Chandra which does not exist in that case.

Our second aim is to answer to a remark made by Varadarajan during the Ramis congress. He
pointed the fact that an invariant eigendistribution, considered as a distribution on the
Lie algebra by the exponential map, is not a solution of the Hotta-Kashiwara module. A
key point in the original proof of Harish-Chandra is precisely the proof that after
multiplication by a function, the eigendistribution is solution of the Hotta-Kashiwara
module (see \cite{VARAD2}). The study of the Hotta-Kashiwara module did not bypass this
difficult step. Here we consider a family a holonomic $\CaD$-module, which we call
(H-C)-modules; this family includes the Hotta-Kashiwara modules but also the module
satisfied directly by an eigendistribution. We prove that these modules are tame and get
a direct proof of Harish-Chandra theorem.

\section{$V$-filtration and $b$-functions.}  \label{chapa1}

We first recall the definition and a few properties of the classical $V$-filtration, then we give a
new definition of quasi-homogeneous $b$-functions and of tame $\CaD$-modules. We end this section
with a result on the inverse image of $\CaD$-modules which will be a key point of the proof in the
next section.

\subsection{Standard $V$-filtrations.}

In this paper, $(X,\OX)$ is a smooth algebraic variety defined over $k$, an algebraically closed
field of characteristic $0$. The sheaf of differential operators with coefficients in $\OX$ is
denoted by $\DX$. The results and the proofs are still valid if $k=\C$, $X$ is a complex analytic
manifold and $\DX$ is the sheaf of differential operators with holomorphic coefficients.

Let $Y$ be a smooth subvariety of $X$ and $\CI_Y$ the ideal of definition of $Y$. The
$V$-filtration along $Y$ is given by \cite{KVAN}:
$$V_k\DX=\ensemble{P\in\DX|_Y}{\forall l\in \Z, P\CI_Y^l\subset \CI_Y^{l+k}}$$
(with $\CI_Y^l=\OX$ if $l\le 0$).

This filtration has been widely used in the theory of $\CaD$-modules, let us recall some
of its properties (for the details, we refer to \cite{SCHAPBOOK}, \cite{ENS},
\cite{SABB}, \cite{IRRMALYL}). The associated graded ring $gr_V\DX$ is the direct image
by $p:T_YX\to X$ of the sheaf $\CaD_{T_YX}$ of differential operators on the normal
bundle $T_YX$. If $\CM$ is a coherent $\DX$-module, a $V\DX$-filtration on $\CM$ is a
good filtration if it is locally finite, i.e. if, locally, there are sections
$(u_1,\dots,u_N)$ of $\CM$ and integers $(k_1,\dots,k_N)$ such that $V_k\CM=\sum
V_{k-k_i}\DX u_i$.

If $\CM$ is a coherent $\DX$-module provided with a good $V$-filtration, the associated graded
module is a coherent $gr_V\DX$-module and if $\CN$ is a coherent submodule of $\CM$ the induced
filtration is a good filtration (see \cite[Chapter III, Proposition 1.4.3]{SCHAPBOOK} or
\cite{SABB}).

Let $\gth_Y$ be the Euler vector field of the fiber bundle $T_YX$, that is the vector
field verifying $\gth_Y(f)=kf$ when $f$ is a function on $T_YX$ homogeneous of degree $k$
in the fibers of $p$. A \textsl{b-function} along $Y$ for a coherent  $\DX$-module with a
good $V$-filtration is a polynomial $b$ such that
$$\forall k\in\Z,\qquad b(\gth_Y+k)gr_V^k\CM=0$$
If the good $V$-filtration is replaced by another, the roots of $b$ are translated by integers.
Here, we always fix the filtration, in particular, if the $\DX$-modules is of the type $\DX/\CI$,
the good filtration will be induced by the canonical filtration of $\DX$.

\subsection{Quasi-homogeneous $V$-filtrations and quasi-$b$-functions.}\label{sec:vfiltr}

Let $\gf=(\gf_1,\dots,\gf_d)$ be a polynomial map from $X$ to the vector space $W=k^d$
and $m_1,\dots,m_d$ be strictly positive and relatively prime integers. We define a
filtration on $\OX$ by:
$$V_k^\gf\OX=\sum_{\scal m\ga=-k}\OX\gf^\ga$$
with $\ga\in\N^d$, $\scal m\ga=\sum m_i\ga_i$ and $\gf^\ga=\gf_1^{\ga_1}\dots\gf_d^{\ga_d}$. If
$k\ge 0$ we set $V_k^\gf\OX=\OX$.

This filtration extends to $\DX$ by:
\begin{equation}\label{equ:vf}
 V_k^\gf\DX=\ensemble{P\in\DX}{\forall l\in \Z, PV_l^\gf\OX\subset V_{l+k}^\gf\OX}
\end{equation}

\begin{defi}
A $(\gf,m)$-weighted Euler vector field is a vector field $\gh$ in $\sum_i \gf_i{\mathcal
V}_X$ such that $\gh(\gf_i)=m_i\gf_i$ for $i=1\dots d$. (${\mathcal V}_X$ is the sheaf of
vector fields on $X$.)
\end{defi}

\begin{lemm}\label{lem:differe}
Any $(\gf,m)$-weighted Euler vector field is in $V_0^\gf\DX$ and if $\gh_1$ and $\gh_2$
are two $(\gf,m)$-weighted Euler vector fields, $\gh_1-\gh_2$ is in $V_{-1}^\gf\DX$.
\end{lemm}

The map $\gf$ may be not defined on $X$ but on an \'etale covering of $X$. More
precisely, let us consider an \'etale morphism $\gn:X'\to X$ and a morphism $\gf:X'\to
W=k^d$. If $m_1,\dots,m_d$ are strictly positive and relatively prime integers, we define
$V_k^\gf\OX$ as the sheaf of functions on $X$ such that $f_\circ\gn$ is in
$V_k^\gf\CO_{X'}$. This defines a $V$-filtration on $\OX$ and on $\DX$ by the formula
(\ref{equ:vf}). The map $TX'\to TX\x_XX'$ is an isomorphism and a vector field $\gh$ on
$X$ defines a unique vector field $\gt^*(\gh)$ on $X'$. So, we define a
$(\gf,m)$-weighted Euler vector field on $X$ as a vector field $\gh$ on $X$ such that
$\gt^*(\gh)$ is a $(\gf,m)$-weighted Euler vector field on $X'$.

\begin{defi}\label{def:section}
Let $u$ be a section of a coherent $\DX$-module $\CM$. A polynomial $b$ is a quasi-$b$-function of
type $(\gf,m)$ for $u$ if there exist a $(\gf,m)$-weighted Euler vector field $\gh$ and a
differential operator $Q$ in $V_{-1}^\gf\DX$ such that $(b(\gh)+Q)u=0$.

The quasi-$b$-function is said \textsl{regular} if the order of $Q$ as a differential operator is
less or equal to the order of the polynomial $b$ and \textsl{monodromic} if $Q=0$.

The quasi-$b$-function is said \textsl{tame} if the roots of $b$ are strictly greater than $-\sum
m_i$.
\end{defi}

These definitions are valid for any map $\gf$ but here we always assume that $\gf$ is smooth. Then
if $Y=\gf^{-1}(0)$, we say for short that $b$ is a quasi-$b$-function of total weight $|m|=\sum
m_i$ along $Y$. Remark that lemma \ref{lem:differe} shows that the definition is independent of the
$(\gf,m)$-weighted Euler vector field $\gh$.

Let $\CM$ be a coherent $\DX$-module. A $V^\gf\DX$-filtration on $\CM$ is a good filtration if it
is locally finite.

\begin{defi}\label{def:filtration}
Let $\CM$ be a coherent $\DX$-module and $V^\gf\CM$ a good $V^\gf\DX$-filtration. A polynomial $b$
is a quasi-$b$-function of type $(\gf,m)$ for $V^\gf\CM$ if, for any $k\in\Z$,
$b(\gh+k)V_k^\gf\CM\subset V_{k-1}^\gf\CM$ where $\gh$ is a $(\gf,m)$-weighted Euler vector field.

The quasi-$b$-function is \textsl{monodromic} if $b(\gh+k)V_k^\gf\CM=0$.
\end{defi}

Definition \ref{def:section} is a special case of definition \ref{def:filtration} if $\DX u$ is
provided with the filtration induced by the canonical filtration of $\DX$.

Recall that if $\CM$ is a $\DX$-module its inverse image by $\gn$ is its inverse image as an
$\OX$-module, that is:
$$\gn^+\CM=\CO_{X'}\otimes_{\gn^{-1}\OX}\gn^{-1}\CM=\DXpX\otimes_{\gn^{-1}\DX}\gn^{-1}\CM$$
where $\DXpX$ is the $(\CaD_{X'},\gn^{-1}\DX)$-bimodule $\CO_{X'}\otimes_{\gn^{-1}\OX}\gn^{-1}\DX$.

\begin{lemm}\label{lem:etale}
Let $\gn:X'\to X$ be an \'etale morphism  and let $\gf$ be a morphism $X'\to W=k^d$. Let
$\CM$ be a coherent $\DX$-module.

The polynomial $b$ is a quasi-$b$-function of type $(\gf,m)$ for a section $u$ of $\CM$ if and only
if it is a quasi-$b$-function of type $(\gf,m)$ for the section $1\otimes u$ of $\gn^+\CM$.
\end{lemm}

\begin{proof}
If $\gn:X'\to X$ is \'etale, the canonical morphism $\DXp\to\DXpX$ given by $P\mapsto P(1\otimes
1)$ is an isomorphism and defines an injective morphism $\gn^*:\gn^{-1}\DX\to\DXp$.

Conversely, the morphism $\tilde\gn:\gn_*\CO_{X'}\to\OX$ given by
$\tilde\gn(f)(x)=\sum_{y\in\gn^{-1}(x)}f(y)$ extends to a morphism $\gn_*\DXp\to\DX$.

These two morphism are compatible with the $V$-filtration defined by $\gf$ and , by definition, a
vector field $\gh$ on $X$ is a $(\gf,m)$-weighted Euler vector field if and only if $\gn^*(\gh)$ is
a $(\gf,m)$-weighted Euler vector field on $X'$. If $(b(\gh)+R)u=0$ we have
$(b(\gn^*\gh)+\gn^*R)(1\otimes u)=0$ and conversely, if $(b(\gn^*\gh)+R_1)(1\otimes u)=0$ then
$(b(\gh)+\gn_*R_1)u=0$.
\end{proof}

\begin{rema}
In \cite{GL} we gave an other definition of the $V^*$-filtration and quasi-$b$-function. The two
definitions are essentially equivalent in the analytic framework but may differ in the algebraic
case. More precisely, the filtration in \cite{GL} is given by a vector field $\gh$ which we called
positive definite. For a given $V^\gf$-filtration, we may find a defining vector field with
coefficients in formal power series (or in convergent series if $k=\C$) but in general not in
rational functions. The definition of \cite{GL} is more intrinsic in the analytic case but not
suitable here.
\end{rema}

\subsection{Tame $\CaD$-modules.}

Let us recall that a stratification of the manifold $X$ is a union $X=\bigcup_\ga X_\ga$ such that

\begin{itemize}
\item For each $\ga$, $\overline{X}_\ga$ is an algebraic subset
  of $X$ and $X_\ga$ is its regular part.

\item $\{X_\ga\}_\ga$ is locally finite.

\item $X_\ga\cap X_\gb =\emptyset$ for$\ga\neq \gb$.

\item If $\overline{X}_\ga\cap X_\gb \neq \emptyset$ then
$\overline{X}_\ga\supset X_\gb$.
\end{itemize}

If $\CM$ is a holonomic $\DX$-module, its characteristic variety $\CHM$ is a homogeneous lagrangian
subvariety of $\TX$ hence there exists a stratification $X=\bigcup X_\ga$ such that
$\CHM\subset\bigcup_\ga \overline{T^*_{X_\ga}X}$ \cite[Ch. 5]{KASHCOUR}. The set of points of $X$
where $\CHM$ is contained in the zero section of $T^*X$ is a non empty Zarisky open subset of $X$,
its complementary is the \textsl{singular support} of $\CM$.

For the next definition, we consider a cyclic $\DX$-module with a canonical generator $\CM=\DX u =
\DX/\CI$ where $\CI$ is a coherent ideal of $\DX$.

\begin{defi}\label{def:quasitame}
The cyclic holonomic $\DX$-module $\CM = \DX u$ is \textsl{tame} if there is a stratification
$X=\bigcup X_\ga$ of $X$ such that $\CHM\subset \bigcup_\ga \overline{T^*_{X_\ga}X}$ and, for each
$\ga$, a tame quasi-$b$-function associated with $X_\ga$.
\end{defi}

With definition \ref{def:section}, this means that for each $\ga$, there is a smooth map $\gf_\ga$
from a Zarisky open set of $X$ to a vector space $V$ such that $X_\ga=\gf_\ga^{-1}(0)$, positive
integers $m_1,\dots,m_d$, a $(\gf,m)$-weighted Euler vector field $\gh$ and a quasi-$b$-function
$b_\ga$ for $u$ with roots $>-\sum m_i$. A subvariety of $X$ is conic for $\gh_\ga$ if it is
invariant under the flow of $\gh_\ga$. The module $\CM$ is \textsl{conic tame} if it satisfy
definition \ref{def:quasitame} and if moreover the singular support of $\CM$ is conic for each
$\gh_\ga$.

The following property of a tame $\DX$-module has been proved in \cite{GL}:

\begin{theo}\label{thm:support}
If the $\DX$-module $\CM$ is tame then it has no quotient with support in a hypersurface
of $X$.
\end{theo}

If $M$ is a real analytic manifold and $X$ its complexification, we also proved:

\begin{theo}\label{thm:lun}
Let $\CM$ be a holonomic and tame $\DX$-module, assume that its singular support is the
complexification of a real subvariety of $M$, then $\CM$ has no distribution solution on
$M$ with support in a hypersurface. If $\CM$ is conic-tame, its distribution solutions
are in $L^1_{loc}$.
\end{theo}

\begin{rema}
It is important to note that the definition of \emph{tame} and the conclusions of theorem
 \ref{thm:lun} depend of the choice of a generator for $\CM$.
\end{rema}

\subsection{Inverse image.}\label{sec:inv}

Let $\gf:X\to W$ and $\gf':X'\to W'$ be two morphisms from smooth algebraic varieties $X$
and $X'$ to the vector spaces $W=k^d$ and $W'=k^{d'}$, let $m_1,\dots,m_d$ and
$m'_1,\dots,m'_{d'}$ be strictly positive integers. Let $f:X'\to X$ and $F:W'\to W$ be
two morphisms such that $\gf_\circ f=F_\circ\gf'$. We assume that $F$ is
\emph{quasi-homogeneous}, that is $F=(F_1,\dots,F_d)$ with
$F_i(\gl^{m'_1}x_1,\dots,\gl^{m'_{d'}}x_{d'})=\gl^{m_i}F(x_1,\dots,x_{d'})$.

If $\CN$ is a $\DX$-module its inverse image by $f$ is:
$$f^+\CN=\CO_{X'}\otimes_{f^{-1}\OX}f^{-1}\CN=\DXpX\otimes_{f^{-1}\DX}f^{-1}\CN$$
where $\DXpX$ is the $(\CaD_{X'},f^{-1}\DX)$-bimodule $\CO_{X'}\otimes_{f^{-1}\OX}f^{-1}\DX$.

We define a filtration on $\DXpX$ by
$$V_k^{\gf'}\DXpX=\sum_{i+j=k}V_i^{\gf'}\CO_{X'}\otimes f^{-1}V_j^\gf\DX$$

By the hypothesis, $g_\circ f$ is a section of $V_k^{\gf'}\CO_{X'}$ for any $g$ section of
$V_k^\gf\OX$, hence the filtration on $\DXpX$ is compatible with the corresponding filtrations on
$\DXp$ and $\DX$.

If a $\DX$-module $\CN$ is provided with a $V^\gf$-filtration, this defines a
$V^{\gf'}\CaD_{X'}$-filtration on $f^+\CN$ by
\begin{equation}\label{equ:filtr}
V_k^{\gf'}f^+\CN=\sum_{i+j=k}V_i^{\gf'}\CO_{X'}\otimes f^{-1}V_j^\gf\CN
=\sum_{i+j=k}V_i^{\gf'}\DXpX\otimes f^{-1}V_j^\gf\CN
\end{equation}

The $V$-filtration has not all the good properties of the usual filtration, in particular non
invertible elements may have an invertible principal symbol. In the proof of theorem
\ref{thm:inverse} we introduce its formal completion given by:
$$\DcXY = \injlim k V_k\DcXY \qquad\textrm{with}
\qquad V_k\DcXY= \projlim lV_k\DX/V_{k-l}\DX$$

By definition the graded ring of $\DcXY$ is the same than the graded ring of $\DX$. If
$\CM$ is a coherent $\DX$-module provided with a good $V$-filtration, its completion
$\widehat V\CM$ is defined in the same way  and has the same associated graded module
than $V\CM$. The following properties may be found in \cite{SCHAPBOOK} and
\cite{IRRMALYL}.

The sheaf $\DcXY$ is a coherent and noetherian, flat over $\DX$. We remind that a
coherent sheaf of rings $\CA$ is noetherian if any increasing sequence of coherent
$\CA$-submodules of a coherent $\CA$-module is stationary. The sheaf of rings $V_0\DcXY$
is also coherent and noetherian.

If $\CM$ is a $\DX$-module provided with a good $V$-filtration, the associated graded
module is a coherent $gr_V\DX$-module and if $\CN$ is a coherent submodule of $\CM$ the
induced filtration is a good filtration. If $\kappa:(\DcXY)^N\to\CM$ is a filtered
morphism which defines a surjective graded morphism $gr_V(\DcXY)^N\to gr_V\CM\to 0$ then
$\kappa$ is surjective.

As $\DcXY$ is flat over $\DX$, if $\CM$ is coherent we have $\widehat
V\CM=\DcXY\otimes_\DX\CM$. Remark also that $\widehat V\OX$, the completion of $\OX$ for
the $V$-filtration, is the formal completion of $\OX$ along $Y$ usually denoted by
$\CO_{\widehat{X|Y}}$ and $\DcXY$ is a $\CO_{\widehat{X|Y}}$-module.

After completion by the $V$-filtration, we get a similar formula:
\begin{equation}\label{equ:filtrf}
\widehat V_k^{\gf'}f^+\CN=\sum_{i+j=k}\widehat V_i^{\gf'}\CO_{X'}\otimes f^{-1}\widehat
V_j^\gf\CN
\end{equation}

 Let $Y=\gf^{-1}(0)$ and $Y'=\gf^{'-1}(0)$, let $p:T_YX\to X$ and $p':T_{Y'}X'\to X'$ be the normal bundles, $\widetilde
f:T_{Y'}X'\to T_YX$ be the map induced by $f$,

\begin{theo}\label{thm:inverse}
We assume that $\gf'$ is smooth on $X'$. If $\CN$ is a holonomic $\DX$-module provided with a good
$V^\gf\DX$-filtration, then $f^+\CN$ is holonomic, ${p'}^{-1}gr_{V^{\gf'}}f^+\CN$ is equal to
$\widetilde f^+ p^{-1}gr_{V^\gf}\CN$ and isomorphic to the graded module associated with a good
$V^{\gf'}\DXp$-filtration of $f^+\CN$.
\end{theo}

\begin{proof}
We recall that if $\CN$ is coherent, then $f^+\CN$ is not coherent in general but if $\CN$ is
holonomic, then $f^+\CN$ is holonomic \cite{INV2}.

The filtration on $\CN$ is a good $V^\gf\DX$-filtration hence we may assume that there are sections
$(u_1,\dots,u_q)$ of $\CN$ and integers $(k_1,\dots,k_q)$ such that $V^\gf_k\CN=\sum
V^\gf_{k-k_i}\DX u_i$. Let $\DXpX[N]$ be the sub-$\DXp$-module of $\DXpX$ generated by the sections
of $\DX$ of order less or equal to $N$. This submodule is finitely generated hence coherent. For
each $N$, $(u_1,\dots,u_q)$ defines a canonical morphism $(\DXpX[N])^q\to f^+\CN$ and the family of
the images of these morphisms is an increasing sequence of coherent submodules of the coherent
$\DXp$-module $f^+\CN$. As $\DXp$ is a noetherian sheaf of rings, this sequence is stationary,
hence there is some $N_0$ such that for each $N>N_0$, the morphism $(\DXpX[N])^q\to f^+\CN$ is
onto. The filtration $V^{\gf'}\DXpX$ induces a good filtration on $\DXpX[N]$ hence, for $N>N_0$ a
good filtration on $f^+\CN$ which is denoted by $V_k^{\gf'}[N]f^+\CN$. To prove the theorem, we
will prove that if $N$ is large enough, $gr_V f^+\CN$ is equal to the graded module $gr_{V[N]}
f^+\CN$ associated with the good filtration $V_k^{\gf'}[N]f^+\CN$.

We assume first that the integers $m'_i$ are equal to $1$, that is that the $V^{\gf'}$-filtration
is the usual $V$-filtration on the non singular variety $Y'=\gf'^{-1}(0)$. For $N>N_0$,
$p'^{-1}gr_{V[N]} f^+\CN$ is a coherent $\CaD_{T_{Y'}X'}$-module. A direct calculation shows that
$p'^{-1}gr_V f^+\CN=\widetilde f^+p^{-1}gr^\gf_V \CN$. If $\CN$ is holonomic then $gr^\gf_V \CN$ is
also holonomic \cite[Cor 4.1.2.]{ENS} hence $p'^{-1}gr_V f^+\CN$ is holonomic hence coherent.

Consider the completion $\widehat V f^+\CN$ of $f^+\CN$ for the $V$-filtration and $\widehat V[N]
f^+\CN$  of $f^+\CN$ for the $V[N]$-filtration. The graded module of $\widehat V f^+\CN$ is equal
to the graded module of $V f^+\CN$ which is coherent. Let $u_1,\dots,u_M$ be local sections of
$\widehat V f^+\CN$ whose classes generate the graded module, then $u_1,\dots,u_M$ generate
$\widehat V f^+\CN$ as a filtered $V\DcXYp$-module and applying the same result to the kernel of
$(V\DcXYp)^M\to \widehat V f^+\CN$ we get that $\widehat V f^+\CN$ admits a filtered presentation
$$(V\DcXYp)^L\to(V\DcXYp)^M\to \widehat V f^+\CN\to 0.$$
This shows in particular that each $\widehat V_k f^+\CN$ is a coherent
$V_0\DcXYp$-module. We know that, for any $N$, $gr_{V[N]} f^+\CN$ is coherent hence for
the same reason, each $\widehat V_k[N] f^+\CN$ is a coherent $V_0\DcXYp$-module.

Consider the family of the images of $\widehat V_k[N] f^+\CN$ in $\widehat V_k f^+\CN$, it is an
increasing sequence of coherent sub-modules of the coherent $V_0\DcXYp$-module $\widehat V_k
f^+\CN$ hence it is stationary because the sheaf of rings $V_0\DcXYp$ is noetherian. Moreover, the
filtration $\widehat V f^+\CN$ is separated hence the maps $\widehat V_k[N] f^+\CN\to \widehat V_k
f^+\CN$ are injective and the union of the images is all $\widehat V_k f^+\CN$, so there is some
$N_0$ such that for any $N>N_0$, $\widehat V_k[N] f^+\CN= \widehat V_k f^+\CN$. This implies that
$gr_V f^+\CN=gr_{V[N]} f^+\CN$ is the graded module associated with a good $V$-filtration of
$f^+\CN$.

Assume now that the numbers $m'_i$ are positive integers. Let $W''=W'$, we define the
ramification map $F_m:W''\to W'$ by $F(s_1,\dots,s_d)=(s_1^{m'_1}, \dots, s_d^{m'_d})$
and the corresponding map $f_m:X''=X'\times_{W'}W''\to X$. Applying the first part of the
proof, we get $\widehat V[N] f_m^+f^+\CN = \widehat V f_m^+f^+\CN$ if $N$ is large. The
formula (\ref{equ:filtrf}) shows that
\begin{equation*}
\widehat V f_m^+f^+\CN=\widehat V\CO_{X''}\otimes_{f^{-1}\widehat V\CO_{X'}}
f^{-1}\widehat V_j^\gf\CN= \CO_{\widehat W}\otimes_{\CO_{\widehat{V'}}} f^{-1}\widehat
V_j^\gf\CN.
\end{equation*}
Here $\CO_{\widehat W}$ is the set of formal power series in $(s_1,\dots,s_d)$ while
$\CO_{\widehat{V'}}$ is the set of formal power series in $(s_1^{m'_1}, \dots,
s_d^{m'_d})$ hence $\CO_{\widehat W}$ is a finite free $\CO_{\widehat{V'}}$-module. So,
if $M'=\sum m'_i$, $\widehat V f_m^+f^+\CN$ is isomorphic to $(\widehat V f^+\CN)^{M'}$
as a $\widehat V\CO_{X'}$-module.

In the same way, $\widehat V[N] f_m^+f^+\CN$ is isomorphic to $(\widehat V[N] f^+\CN)^{M'}$, hence
$\widehat V[N] f^+\CN=\widehat V[N] f^+\CN$. This shows that $gr_{V^{\gf'}}f^+\CN$ is the graded
module associated with $V^{\gf'}[N] f^+\CN$ which is a good filtration of $f^+\CN$.
\end{proof}

\begin{rema}
The result was known when $f$ is a submersion, $Y'=f^{-1}(Y)$ and the $V$-filtrations
being the usual $V$-filtrations along $Y$ and $Y'$ \cite{IRRMALYL}. The introduction of
the weights $m_i$ and $m'_i$ allows $f$ to be non submersive and $Y'$ to be a proper
subvariety of $f^{-1}(Y)$; the relation between the weights is given by the
quasi-homogeneity of $F$.
\end{rema}

\begin{coro} \label{cor:invbf}
Under the hypothesis of theorem \ref{thm:inverse}, if $\CN$ is a holonomic $\DX$-module provided
with a good $V^\gf\DX$-filtration, $f^+\CN$ is provided with a good $V^{\gf'}\DXp$-filtration such
that a polynomial $b$ is a quasi-$b$-function of type $(\gf,m)$ for the filtration of $\CN$ if and
only if $b$ is a quasi-$b$-function of type $(\gf',m')$ for the filtration of $f^+\CN$.
\end{coro}

\begin{proof}
Let $\gh'$ be a $(\gf',m')$-weighted Euler vector field, then $\gh=f_*\gh'$ is a $(\gf,m)$-weighted
Euler vector field. As definition \ref{def:filtration} is independent of the $(\gf,m)$-weighted
Euler vector field, we may assume that the quasi-$b$-function for $\CN$ is relative to $\gh$.

By definition, for any $Q$ in $\DXpX$, we have $\gh'Q = Q\gh$ hence for any polynomial $b(\gh')Q =
Qb(\gh)$ which shows the corollary.
\end{proof}

\begin{coro} \label{cor:invbf2}
Under the hypothesis of theorem \ref{thm:inverse}, if $\CN$ is a holonomic $\DX$-module and $u$ a
section of $\CN$ with a  quasi-$b$-function of type $(\gf,m)$, then  the section $1_{X'\to
X}\otimes u$ of $f^+\CN$ has the same polynomial $b$ as a quasi-$b$-function of type $(\gf',m')$.
\end{coro}

\begin{proof} Recall that $1_{X'\to X}$ is the canonical section $1\otimes 1$ in
$\DXpX=\CO_{X'}\otimes_{f^{-1}\OX}f^{-1}\DX$. If $u$ is a section of $\CN$, we set on $\DX u$ the
filtration image of the filtration of $\DX$. Then, by definition of the filtration
$V^{\gf'}[N]f^+\CN$ used in the proof of theorem \ref{thm:inverse}, $1_{X'\to X}\otimes u$ is of
order $0$ for this filtration. Then corollary \ref{cor:invbf2} is a special case of corollary
\ref{cor:invbf}.
\end{proof}

\section{Reductive Lie algebras.}  \label{chaph}

\subsection{Statement of the main theorem.}\label{sec:HC}

Let $G$ be a connected reductive algebraic group with Lie algebra $\kg$, let $\kg^*$ be the dual
space of $\kg$. The group $G$ acts on $\kg$ by the adjoint action hence on the symmetric algebra
$S(\kg)$ identified to the space $\CO(\kg^*)$ of polynomial functions on $\kg^*$. By Chevalley's
theorem, the space $\CO(\kg^*)^G\simeq S(\kg)^G$ of invariant polynomials on $\kg^*$ is equal to a
polynomial algebra $k[Q_1,\dots,Q_l]$ where $Q_1,\dots,Q_l$ are algebraically independent. These
spaces are graded and we denote $S_+(\kg)^G=\oplus_{k>0}S_k(\kg)^G$. It is also the set
$\CO_+(\kg^*)^G$ of invariant polynomials vanishing at $\{0\}$, their common roots are the
nilpotent elements of $\kg^*$.

The differential of the adjoint action induces a Lie algebra morphism $\gt:\kg\to
\mathrm{Der}S(\kg^*)$ by:
$$(\gt(A)f)(x)=\frac d{dt}f\left(\exp(-tA).x\right)|_{t=0}
\quad\mathrm{for}\quad A\in\kg, f\in S(\kg^*)=\CO(\kg), x\in\kg$$ i.e. $\gt(A)$ is the vector field
on $\kg$ whose value at $x\in\kg$ is to $[x,A]$. We denote by $\gt(\kg)$ the set of all vector
fields $\gt(A)$ for $A\in\kg$. It is the set of vector fields on $\kg$ tangent to the orbits of
$G$.

Let $\Dgg$ be the sheaf of differential operators on $\kg$ invariant under the adjoint action of
$G$. The principal symbol $\gs(P)$ of such an operator $P$ is a function on $T^*\kg=\kg\x\kg^*$
invariant under the action of $G$. Examples of such invariant functions are the elements of
$S(\kg)^G$ identified to functions on $\kg\x\kg^*$ constant in the variables of $\kg$. If $F$ is a
subsheaf of $\Dgg$, we denote by $\gs(F)$ the sheaf of the principal symbols of all elements of
$F$.

\begin{defi}
A subsheaf $F$ of $\Dgg$ is of (H-C)-type if $\gs(F)$ contains a power of $S_+(\kg)^G$.
An (H-C)-type $\Dg$-module is the quotient $\MF$ of $\Dg$ by the ideal $\CI_F$ generated
by $\gt(\kg)$ and  by $F$.
\end{defi}

The main result of this paper is

\begin{theo}\label{thm:main}
Any $\Dg$-module of (H-C)-type is holonomic and conic-tame.
\end{theo}
Here (H-C) stands for Harish-Chandra. There are two main examples of such $\Dg$-modules which we
describe now.
\begin{exem}
An element $A$ of $\kg$ defines a vector field with constant coefficients on $\kg$ by:
$$(A(D_x)f)(x)=\frac d{dt}f(x+tA)|_{t=0}\quad\mathrm{for}\quad f\in S(\kg^*), x\in\kg$$

By multiplication, this extends to an injective morphism from $S(\kg)$ to the algebra of
differential operators with constant coefficients on $\kg$; we identify $S(\kg)$ with its image and
denote by $P(D_x)$ the image of $P\in S(\kg)$. If $F$ is a finite codimensional ideal of
$S(\kg)^G$, its graded ideal contains a power of $S_+(\kg)^G$ hence when it is identified to a set
of differential operators with constant coefficients, $F$ is a subsheaf of $\Dg$ of (H-C)-type and
$\MF$ is a $\Dg$-module of (H-C)-type. If $\gl\in\kg^*$, the module $\MFl$ defined by Hotta and
Kashiwara \cite{HOTTA} is the special case where $F$ is the set of polynomials $Q-Q(\gl)$ for $Q\in
S(\kg)^G$.
\end{exem}

\begin{exem}
The enveloping algebra $U(\kg)$ is the algebra of left invariant differential operators
on $G$. It is filtered by the order of operators and the associated graded algebra is
isomorphic by the symbol map to $S(\kg)$. This map is a $G$-map and defines a morphism
from the space of bi-invariant operators on $G$ to the space $S(\kg)^G$. This map is a
linear isomorphism, its inverse is given by a symmetrization morphism \cite[Theorem
3.3.4.]{VARADA}. We assume that $k=\C$. Then, through the exponentional map a
bi-invariant operator $P$ defines a differential operator $\widetilde P$ on the Lie
algebra $\kg$ which is invariant under the adjoint action of $G$ (because the exponential
intertwines the adjoint action on the group and on the algebra) and the principal symbol
$\gs(\widetilde P)$ is equal to $\gs(P)$.

Let $U$ be an open subset of $\kg$ where the exponential is injective and $U_G=\exp(U)$. Let $T$ be
an invariant eigendistribution on $U_G$ and $\widetilde T$ the distribution on $U$ given by $\scal
T \gf = \scal{\widetilde T}{\gf_o\exp}$. As $T$ is invariant and eigenvalue of all bi-invariant
operators, $\widetilde T$ is solution of an (H-C)-type $\Dg$-module. As this module is conic-tame
by theorem \ref{thm:main}, the results of theorems \ref{thm:support} and \ref{thm:lun} are true for
it, hence $\widetilde T$ and $T$ are a $L^1_{loc}$-function.
\end{exem}

As $\kg$ is reductive, it is the direct sum $\kc\oplus [\kg,\kg]$ of its center and of the
semi-simple Lie algebra $[\kg,\kg]$. We choose a non-degenerate $G$-invariant symmetric bilinear
form $\gk$ on $\kg$ which extend the Killing form of $[\kg,\kg]$. This defines an isomorphism from
$\kg^*$ to $\kg$ and the cotangent bundle $T^*\kg=\kg\x\kg^*$ is identified with $\kg\x\kg$. Then
if $\CN(\kg)$ is the nilpotent cone of $\kg$, the characteristic variety of an (H-C)-type
$\Dg$-module is a subset of:
$$\ensemble{(x,y)\in\kg\x\kg}{[x,y]=0, y\in \CN(\kg)}$$
so it is a holonomic $\Dg$-module \cite{HOTTA}.

\subsection{Stratification of a reductive Lie algebra.}\label{sec:stratification}

In this section, we define the stratification which will be used to prove that an (H-C)-type module
is tame. This stratification is classical (see \cite{ATI} for example).

The stratification of a reductive Lie algebra is the direct sum of the center by the stratification
of the semi-simple part, so we may assume that $\kg$ is semi-simple. An element $X$ of $\kg$ is
said to be \textsl{semisimple} (resp. \textsl{nilpotent}) if $\ad(X)$ is semisimple (resp $\ad(X)$
is nilpotent). Any $X\in\kg$ may be decomposed in a unique way as $X=S+N$ where $S$ is semisimple,
$N$ is nilpotent and $[S,N]=0$ (Jordan decomposition). An element $X$ is said to be
\textsl{regular} if the dimension of its centralizer $\kg^X=\ensemble{Z\in\kg}{[X,Z]=0}$ is
minimal, that is equal to the rank of $\kg$. The set $\kgrs$ of semisimple regular elements of
$\kg$ is Zarisky dense and its complementary $\kg'$ is defined by a $G$-invariant polynomial
equation $\gD(X)=0$. The function $\gD$ may be defined from the characteristic polynomial of
$\ad(X)$:
$$\det(T.Id-ad(X))=T^n+\sum \gl_i(X)T^i$$
Here $n$ is the dimension of $\kg$. Then $\gl_0\equiv 0$, the rank $l$ of $\kg$ is the
lowest $i$ such that $\gl_i\not\equiv 0$ and $\gD(X)=\gl_l(X)$. This function is
homogeneous of degree $n-l$.

The set $\kN(\kg)$ of nilpotent elements of $\kg$ is a cone equal to:
$$\kN(\kg)=\ensemble{X\in\kg}{\forall P\in\CO(\kg)^G\quad P(X)=P(0)}$$
and the set of nilpotent orbits is finite and define a stratification of $\kN$ \cite[Cor
3.7.]{KOS1}.

We fix a Cartan subalgebra $\kh$ of $\kg$ and denote by $\kw$ the Weyl group $\kw(\kg,\kh)$. Let
$\gF=\gF(\kg,\kh)$ be the root system associated with $\kh$. For each $\ga\in\gF$ we denote by
$\kg_\ga$ the root subspace corresponding to $\ga$ and by $\kh_\ga$ the subset
$[\kg_\ga,\kg_{-\ga}]$ of $\kh$ (they are all $1$-dimensional). Let $\CF$ be the set of the subsets
$P$ of $\gF$ which are closed and symmetric that is such that $(P+P)\cap\gF\subset P$ and $P=-P$.
For each $P\in\CF$ we define $\khp=\sum_{\ga\in P}\kh_\ga$, $\kgp=\sum_{\ga\in P}\kg_\ga$,
$\khpo=\ensemble{H\in\kh}{\ga(H)=0 \textrm{ if } \ga\in P}$ and $\khpp =\ensemble{H\in\kh}{\ga(H)=0
\textrm{ if } \ga\in P,\ga(H)\ne0 \textrm{ if } \ga\notin P}$.

The following results are well-known (see \cite[Ch VIII, \S3]{BOU}):

a) $\kqp=\khp+\kgp$ is a semisimple Lie subalgebra of $\kg$ stable under $\ad\kh$ and
$\khpo$ is an orthocomplement of $\khp$ for the Killing form, $\khp$ is a Cartan
subalgebra of $\kqp$. The Weyl group $\kw_P$ of $(\kqp,\khp)$ is identified to the
subgroup $\kw'$ of $\kw$ of elements whose restriction to $\khpo$ is the identity.

b) $\kh+\kgp$ is a reductive Lie subalgebra of $\kg$ stable under $\ad\kh$. For any
$S\in\khpo$, $\kh+\kgp\subset\kg^S$ and $\khpp=\ensemble{S\in\khpo}{\kg^S=\kh+\kgp}$.

c) Conversely, if $S\in\kh$, there exists a subset $P$ of $\gF$ which is closed and
symmetric such that $\kg^S=\kh+\kgp$. $P$ is unique up to a conjugation by $\kw$.

To each $P$ of $\CF$ and each nilpotent orbit $\kO$ of $\kqp$ we associate a conic subset
of $\kg$
\begin{equation}
S_{(P,\kO)}=\bigcup_{X\in\khpp}G.(X+\kO)\label{def:strat}
\end{equation}
where $G.(X+\kO)$ is the union of orbits of points $X+\kO$.

If $X=S+N$ is the Jordan decomposition of $X\in\kg$, the semisimple part $S$ belongs to a
Cartan subalgebra which we may assume to be $\kh$ because they are all conjugate. Hence
there is some $P$ in $\CF$ such that $\kg^S=\kh+\kgp$. Then, if the orbit of $N$ in
$\kqp$ is $\kO$, $X\in S_{(P,\kO)}$. For a detailed proof of the fact that it is a
stratification, see \cite{GL} .

\subsection{Polynomials and differentials.}\label{sec:inver}

Let us begin with some elementary calculations. If
$\gb=(\gb_1,\dots,\gb_n)$ is a multi-index of $\N^n$ we denote
$|\gb|=\sum\gb_i$ and $\gb!=\gb_1!\dots\gb_n!$, if $\ga$ is another
element of $\N^n$, we denote by $\ga\le\gb$ the relation
$\ga_1\le\gb_1,\dots,\ga_n\le\gb_n$.

\begin{lemm}\label{lem:binom}
Let $\gb\in\N^n$ and $M=|\gb|$, let $N\in\N$ such that $N\le M$, then
$$\sum_{\substack{|\ga|=N\\ \ga\le\gb}}
 \frac{\gb!}{\ga!(\gb-\ga)!}=\frac{M!}{N!(M-N)!}$$
\end{lemm}

\begin{proof}
$$\sum_{\ga\le\gb} \frac{\gb!}{\ga!(\gb-\ga)!}x^\ga=
\prod_{i=1}^n\sum_{\ga_i=0}^{\gb_i}\frac{\gb_i!}{\ga_i!(\gb_i-\ga_i)!}x_i^{\ga_i}
=(1+x_1)^{\gb_1}\dots(1+x_n)^{\gb_n}$$
hence if $t=x_1=\dots=x_n$ we get~:
$$\sum_{\ga\le\gb} \frac{\gb!}{\ga!(\gb-\ga)!}t^{|\ga|}=(1+t)^M$$
and the coefficient of $t^N$ in both side of the equality gives the lemma.
\end{proof}

\begin{lemm}\label{lem:bino2}
Let us denote $x=(x_1,\dots,x_n)$, $D_x=(D_{x_1},\dots,D_{x_n})$,
$x^\ga=(x_1^{\ga_1},\dots,x_n^{\ga_n})$ and $D_x^\ga=(D_{x_1}^{\ga_1},\dots,D_{x_n}^{\ga_n})$, let
$\gth=\sum x_iD_{x_i}$, then~:
$$\sum_{|\ga|=N}\frac{N!}{\ga!}x^\ga D_x^\ga=\gth(\gth-1)\dots(\gth-N+1)$$
\end{lemm}

\begin{proof}
To prove the equality of the two differential operators we have to show that they give the same
result when acting on a monomial $x^\gb$, so lemma \ref{lem:binom} gives~:
$$\sum_{|\ga|=N}\frac{N!}{\ga!}x^\ga D_x^\ga x^\gb=
\sum_{\substack{|\ga|=N\\ \ga\le\gb}}\frac{N!}{\ga!}\frac{\gb!}{(\gb-\ga)!} x^\gb
=\frac{|\gb|!}{(|\gb|-N)!} x^\gb=\gth(\gth-1)\dots(\gth-N+1)x^\gb$$
\end{proof}

\begin{prop}\label{prop:bff}
Let $p_1,\dots,p_n(\gx)$ be homogeneous polynomial on $X=\C^n$ and assume that:
$$\bigcap_{i=1}^n\{p_i(\gx)=0\}=\{0\}$$
Let $\CI$ be the ideal of $\DX$ generated by $p_1(D_x)\dots p_n(D_x)$ and $\CM=\DX/\CI$. The $\DX$
module $\CM$ is holonomic and the $b$-function of $\CM$ relative to $\{0\}$ is equal to
$$b(\gth)=\gth(\gth-1)\dots(\gth+n-M)$$
where $M$ is the sum of the degrees of the polynomials $p_1,\dots,p_n$ and $\gth$ the
Euler vector field of $X$. This $b$-function is monodromic in the canonical coordinates
of $\C^n$.
\end{prop}

\begin{proof}
The Nullstellensatz shows that there is some integer $M_1$ such that the monomial
$\gx^\ga$ are in the ideal generated by $p_1,\dots,p_n$ if $\abs\ga>M_1$. In fact it is
known that the lowest $M_1$ is $M-n$ (the proof uses the Hilbert polynomial). Then lemma
\ref{lem:bino2} shows that the $b$-function of $\CM$ divides $\gth\dots(\gth+n-M)$. It
has been proved by T.Torrelli \cite{TORR} that all integers $0,\dots,M-n$ appear
effectively as roots of $b$.
\end{proof}

\begin{prop}\label{prop:bf}
Let $p_1,\dots,p_n$ be the same polynomials as in the previous proposition and let $P_1,\dots P_n$
be differential operators such that $\gs(P_i)=p_i$. Let $\CI$ be the ideal of $\DX$ generated by
the operators $P_1,\dots,P_n$ and $\CM=\DX/\CI$. The $\DX$-module $\CM$ is holonomic and the
$b$-function of $\CM$ relative to $\{0\}$ is equal to
$$b(\gth)=\gth(\gth-1)\dots(\gth+n-M)$$
The $b$-function of $\CM$ along a vector subspace $L$ of $\C^n$ divides the same polynomial $b$.
\end{prop}
\begin{proof}
Each function $\gx^\ga$ for $|\ga|=N=M-n+1$ is written as $\gx^\ga=\sum q_i^\ga(\gx)p_i(\xi)$ and
\begin{align*}
b(\gth)=&\sum_{|\ga|=N}\frac{N!}{\ga!}x^\ga D_x^\ga=
\sum_{\substack{|\ga|=N\\ i=1,\dots,n}}\frac{N!}{\ga!}x^\ga
q_i^\ga(D_x)p_i(D_x)\\
=&\sum_{\substack{|\ga|=N\\i=1,\dots,n}}\frac{N!}{\ga!}x^\ga q_i^\ga(D_x)P_i(x,D_x)
+\sum_{\substack{|\ga|=N\\ i=1,\dots,n}}\frac{N!}{\ga!}x^\ga
q_i^\ga(D_x)\left(p_i(D_x)-P_i(x,D_x)\right)
\end{align*}
By definition, the order for the $V$-filtration along $\{0\}$ is always less then the usual order
with equality if the operator has constant coefficients. So the order for the $V$-filtration of
$\sum\frac{N!}{\ga!}x^\ga q_i^\ga(D_x)(p_i(D_x)-P_i(x,D_x))$  is strictly less than the order of
$\sum\frac{N!}{\ga!}x^\ga q_i^\ga(D_x)p_i(D_x)$ which is the order of $b(\gth)$, that is $0$. On
the other hand $\sum\frac{N!}{\ga!}x^\ga q_i^\ga(D_x)P_i(x,D_x)$ is in the ideal $\CI$ hence
$b(\gth)$ is a $b$-function.

For the second part, we choose linear coordinates of $\C^n$ such
that $L=\ensemble{(x,t_1,\dots,t_d)\in\C^n}{t=0}$ and we write
\[b(<t,D_t>)=\sum_{\abs\gb=N}\frac{N!}{\gb!}t^\gb D_t^\gb\]
As all $D_t^\gb$ for $\abs\gb=N$ are in the ideal generated by
$p_1(D_x,Dt),\dots,p_n(D_x,Dt)$ the proof is the same then before.
\end{proof}

\subsection{Proof of the main theorem.}\label{sec:proof}

Let $\gf:Y\to X$ be an algebraic map. A vector field $u$ on $Y$ is said to be \textsl{tangent to
the fibers of $\gf$} if $u(f\circ\gf)=0$ for all $f$ in $\OX$. A differential operator $P$ is said
to be \textsl{invariant under $\gf$} if there exists a $k$-endomorphism $A$ of $\OX$ such that
$P(f\circ\gf)=A(f)\circ\gf$ for all $f$ in $\OX$. If we assume from now that $\gf$ is dominant, $A$
is uniquely determined by $P$ and is a differential operator on $X$. We denote by $A=\gf_*(P)$ the
image of $P$ in $\DX$ under this ring homomorphism.

We fix a Cartan subalgebra $\kh$ of $\kg$ and denote by $\kw$ the Weyl group
$\kw(\kg,\kh)$. The Chevalley theorem shows that $\CO(\kg)^G$ is equal to
$k[P_1,\dots,P_l]$ where $(P_1,\dots,P_l)$ are algebraically independent invariant
polynomials and $l$ is the rank of $\kg$, that the set of polynomials on $\kh$ invariant
under $\kw$ is $\CO(\kh)^\kw=k[p_1,\dots,p_l]$ where $p_j$ is the restriction to $\kh$ of
$P_j$ and that the restriction map $P\mapsto P|_\kh$ defines an isomorphism of
$\CO(\kg)^G$ onto $\CO(\kh)^\kw$ \cite[\S 4.9.]{VARADA}. The space $W=\kh/\kw$ is thus
isomorphic to $k^l$ and the functions $P_1,\dots,P_l$ define two morphisms $\gy:\kg\to W$
and $\gf:\kh\to W$ by $\gy(x)=(P_1(x),\dots,P_l(x))$ and $\gf(z)=(p_1(z),\dots,p_l(z))$.

An operator $Q$ of $\Dg^G$ transforms invariant functions into invariant functions hence
is invariant under $\gy$ and $\gy_*(Q)$ is a differential operator on $W$. A vector field
of $\gt(\kg)$ annihilates the functions $P_1,\dots,P_l$ hence is tangent to the fibers of
$\gy$. In the same way, let $\Dh^\kw$ be the space of differential operators on $\kh$
which are invariant under the action of the Weyl group $\kw$, they are invariant under
$\gf$ and define operators on $W$ through $\gf_*$.

Let $\CaD[(\kg)^G$ (resp. $\CaD(\kh)^\kw$) be the set of global sections of $\Dgg$ (resp. of
$\Dh^\kw$). The morphism of Harish-Chandra \cite{HC2} is a morphism of sheaves of rings
$\gd:\CaD(\kg)^G\to\CaD(\kh)^\kw$ which satisfies the following properties:
\begin{enumerate}
    \item If $f\in\CO(\kg)^G\simeq\CO(\kh)^\kw$ then
$\gd(P)(f|_\kh)=\gD^{1/2}P(f)\gD^{-1/2}|_\kh$.
    \item If $f\in\CO(\kg)^G$, $\gd(f)$ is the restriction of $f$ to $\kh$
    \item If $f\in S(\kg)^G$ and $f$ is considered as a constant coefficients operator, then
 $\gd(f)$ is the restriction of $f$ to $\kh^*$.
    \item The morphism $\gd$ is surjective onto $\CaD(\kh)^\kw$.
    \item The kernel of $\gd$ is $\CaD(\kg)^G\cap\CaD(\kg)\gt(\kg)$.
\end{enumerate}

The last two results have been proved algebraically by Levasseur and Stafford in
\cite{LEVASS} and \cite{LEVASS3}. Let $E$ be the Euler vector field of $\kg$ and $\gvt$
the Euler vector field of $\kh$. The function $\gD$ is homogeneous of degree $n-l$
(\ref{sec:stratification}) hence $\gd(E)$ is equal to $\gvt-(n-l)/2$.

Let $\DV[d^{-1}]$ be the sheaf of differential operators on $W$ with poles on $\{d=0\}$ and
$\CaD(W)[d^{-1}]$ be the ring of its global sections. The function $\gD$ is invariant hence of the
form $d(P_1,\dots,P_l)$ and the formula $Q\mapsto d^{1/2}Q d^{-1/2}$ defines an isomorphism $\gg$
of $\DV[d^{-1}]$. We get a diagram:

\begin{equation} \label{HCcomm}
\begin{CD}
\CaD(\kg)^G        @>{\gd}>> \CaD(\kh)^\kw \\
@V{\gy_*}VV                 @V{\gf_*}VV \\
\CaD(W)[d^{-1}]   @>\gg>> \CaD(W)[d^{-1}]\\
\end{CD}
\end{equation}

If $f$ is a polynomial on $W$ and $Q$ an operator of $\Dg^G$ we have
$\gf_*(\gd(P))(f)=\gg(\gy_*(P))(f)$ from the definitions hence the diagram is
commutative. We can avoid the denominators $[d^{-1}]$ in the diagram because of the
following lemma:

\begin{lemm}\label{lem:pasden}
The morphism $\gg$ sends the image of $\gy_*$ into $\CaD(W)$ while its inverse $\gg^{-1}$ sends the
image of $\gf_*$ into $\CaD(W)$.
\end{lemm}

\begin{proof}
This commutativity of the diagram shows that if an operator of $\CaD(W)$ is in the range of $\gy_*$
then its image under $\gg$ is in $\CaD(W)$.

Conversely let us choose a positive system of roots for $(\kg,\kh)$ and define a function
by $\gp=\Pi_{\ga>0}\ga$. Then $\gp$ is a product of distinct linear forms, its square
$\gp^2$ is equal to the restriction of $\gD$ to $\kh$ and it is changed to $-\gp$ under a
reflection of the Weyl group.

Let $P\in\CaD(\kh)^\kw$ and $f\in\CO_\kh^\kw$, by definition the function $Pf$ is invariant under
$\kw$ while the function $\gp^{-1}P(\gp f)$ is in $\CO_\kh[\gp^{-1}]$ and is invariant under $\kw$.
Hence the function $\gt=P(\gp f)$ is in $\CO_\kh$ and changes its sign under the action of
reflections.

Let $z$ be a point of $\{\gp=0\}$, there exists a root $\ga$ such that $\ga(z)=0$. Let
$s$ be the reflection which let the hyperplane $\{\ga=0\}$ invariant. We have
$\gt(z)=\gt(z^s)=-\gt(z)$ hence $\gt(z)=0$. As $\gt$ vanishes on $\{\gp=0\}$ and $\gp$
has multiplicity $1$, $\gt$ is divisible by  $\gp$ and $\gp^{-1}P(\gp f)$ has no
denominator.

So the operator $\gp^{-1}P\gp$ is in $\CaD(\kh)^\kw[\gp^{-1}]$ but applied to an invariant
polynomial it gives a polynomial. Its image under $\gf_*$ is thus a differential operator of
$\CaD(W)[d^{-1}]$ which sends any polynomial to a polynomial hence an operator of $\CaD(W)$.
\end{proof}

Let $F$ be an (H-C)-type subsheaf of $\Dg^G$, we define four $\CaD$-modules:

\begin{itemize}
\item $\MF$ is the (H-C)-type $\Dg$-module. It is equal to the quotient of $\Dg$ by the ideal $\IF$
generated by $\gt(\kg)$ and $F$.
\item $\Nl$ is the quotient of $\DV$ by the ideal generated by $\gy_*(F)$.
\item $\Mlh$ is the quotient of $\Dh$ by the ideal generated by $\gd(F)$.
\item $\Nlh$ is the quotient of $\DV$ by the ideal generated by  $\gf_*(\gd(F))$.
\end{itemize}

Let $1_{\kg\to W}$ be the canonical generator of $\DGV$ as defined in the proof of
corollary \ref{cor:invbf2} and $u_{\kg\to W}$ its class in $\gy^+\Nl$. We denote by $\NO$
the $\Dg$-submodule
 of $\gy^+\Nl$ generated by $u_{\kg\to W}$.
\begin{theo}\label{thm:second}
The module $\NO$ is conic-tame.
\end{theo}

In this section we prove this theorem and in the next section we prove that $\MF$ is isomorphic to
$\NO$.

\begin{prop}\label{prop:Euler}
Let $n$ be the dimension of $\kg$, $l$ its rank. Then there exit some positive integer $N$ such
that $$ b(T)=(T-N)\dots T(T+1)\dots(T+\frac{n-l}2)$$ is a quasi-$b$-function of total weight
$(n+l)/2$ for $\Nl$ along $\{0\}$. Moreover, $N=0$ if $\gs(F)=S_+(\kg)$.
\end{prop}

\begin{proof}We recall that the rank $l$ of the algebra $\kg$ is the dimension of a Cartan subalgebra and that
the degrees $n_1,\dots,n_l$ of the generators $P_1,\dots,P_l$ of $\CO(\kg)^G$ are called
the primitive degrees of $\kg$ and that their sum is $\frac{n+l}2$ \cite{VARADA}. The map
$\gy:\kg\to W$ is defined by $(P_1,\dots,P_l)$, hence if $E=\sum x_iD_{x_i}$ is the Euler
vector field of $\kg$, $\gh=\gy_*(E)$ is equal to $\sum n_it_iD_{t_i}$.

The morphism $\gd$ is graded and its restriction to $S(\kg)^G$ is the map $Q\mapsto
q=Q|_\kh$ hence $\gs(\gd(F))$ the set of principal symbols contains a power of
$S_+(\kh)^\kw$ (and is equal to $S_+(\kh)^\kw$ if  $\gs(F)=S_+(\kg)$). We may then apply
proposition \ref{prop:bf} to the module $\Mlh$ and we find that its $b$-function is equal
to $b_0(\gvt)=\gvt(\gvt-1)\dots(\gvt-M)$ where $\gvt$ is the Euler vector field of $\kh$
and $M$ is a positive integer equal to $(n-l)/2$ if $\gs(\gd(F))=S_+(\kh)^\kw$. This
means that there exist differential operators $R,A_1,\dots,A_l$ on $\kh$ such that $R$ is
of order $-1$ for the $V$-filtration in $\{0\}$ such that
\[b_0(\gvt)+R(z,D_z)=A_1(z,D_z)q_1(z,D_z)+\dots+A_l(z,D_z)q_l(z,D_z)\]
The action of $\kw$ on $\Dh$ does not affect the $V$-filtration and  $b_0(\gvt)$ and all
$q_i(z,D_z)$ are invariant under the Weyl group, so if we take the mean value (that is
$\frac1{\#\kw}\sum_{w\in \kw}P^w$) we find the same relation with $R$ and all $A_i$
invariant under $\kw$.

Applying $\gf_*$ and $\gg^{-1}$ we find
\begin{equation}\label{equ:b0}
b_0(\gg^{-1}(\gf_*(\gvt)))+\gg^{-1}(\gf_*(R))=B_1\gy_*(Q_1)+\dots+B_l\gy_*(Q_l)
\end{equation}
with $B_1,\dots,B_l$ in $\CaD(W)$ (lemma \ref{lem:pasden}) and
$\gg^{-1}(\gf_*(q_i))=\gy_*(Q_i)$ (the diagram \ref{HCcomm} is commutative).

As $\gf=(p_1,\dots,p_l)$ and $p_i$ has degree $n_i$, $\gf_*(\gvt)$ is equal to $\gh=\sum
n_it_iD_{t_i}$. We have $\gg^{-1}(\gh)=d^{-1/2}\gf_*(\gvt)d^{1/2}=\gf_*(\gD^{-1/2}\gvt\gD^{1/2})$
and the function $\gD$ is homogeneous of degree $n-l$ hence $\gD^{-1/2}\gvt\gD^{1/2}=\gvt+(n-l)/2$
and $\gg^{-1}(\gh)=\gh+(n-l)/2$.
\end{proof}

\begin{prop}\label{prop:Euler2}
For each nilpotent orbit $S$ of codimension $r$, $\NO$ has a $b$-function of total weight $(n+r)/2$
along $S$ equal to
$$b(T)=(T-N)\dots T(T+1)\dots(T+\frac{n-l}2)$$
with $N=0$ if $\gs(F)=S_+(\kg)$. Here $n$ is the dimension of $\kg$ and $l$ its rank.

All roots of $b$ are strictly greater than $-\frac{n+r}2$ hence this $b$-function is
tame.
\end{prop}

\begin{proof}
Let us consider first the null orbit $S=\{0\}$. We apply corollary \ref{cor:invbf} to
$X'=\kg$, $X=W$, $f=\gy$, $W'=\kg$, $\gf'$ is the identity map of $\kg$, $\gf$ the
identity map of $W$ and $F:\kg^Y\to W$ is the map $\gy$. The weights $(m'_1,\dots,m'_n)$
on $\kg$ are $(1,\dots,1)$, that is the $V$-filtration on $\kg$ is the usual
$v$-filtration relative to $\{0\}$, and the weights $(m_1,\dots,m_l)$ on $W$ are the
primitive degrees considered in the proof of Proposition \ref{prop:Euler}. Then $F$ is
quasi-homogeneous and we get directly the result for $S=\{0\}$.

Consider now the nilpotent orbit $S$ of maximal dimension, then $S$ is the smooth part of
the nilpotent cone and $\gy:\kg\to W$ is smooth on $S$. We apply corollary
\ref{cor:invbf} to $X'=\kg$, $X=W$, $f=\gy$, $W'=W$, $\gf'=\gy$, $\gf$ and $F$ are both
the identity map of $W$. The weights on $\kg$ and on $W$ are the weights
$(m_1,\dots,m_l)$ considered on $W$ in the case of the null orbit.

We consider now a non null nilpotent orbit $S$. Let $X\in S$, by the Jacobson-Morozov
theorem, we can find $H$ and $Y$ in $\kg$ such that $(H,X,Y)$ is a $\sld$-triple. They
generate a Lie algebra isomorphic to $\sld$ which acts on $\kg$ by the adjoint
representation. The theory of $\sld$-representations shows $\kg$ splits into a direct sum
$\bigoplus_{i=1}^r E(\gl_i)$ of irreducible submodules. The dimension of $E(\gl_i)$ is
$\gl_i+1$ hence $n=\sum (\gl_i+1)$. Moreover $\kg=[X,\kg]\oplus\kg^Y$, $\dim \kg^Y = r$
and we can select a basis $(Y_1,\dots,Y_r)$ of $\kg^Y$ such that $[H,Y_i]=-\gl_iY_i$. The
tangent space to $S$ at $X$ is $[X,\kg]$ hence $r$ is the codimension of $S$.

The map $\gn:G\x \kg^Y\to \kg$ given by $\gn(g,Z)=g.(X+Z)$ is a submersion because its
tangent map is the map $\kg\x\kg^Y\to\kg$ given by $(Z',Z)\mapsto [Z',X]+Z$. Let $\kg_1$
be a linear subspace of $\kg$ such that $\kg=\kg^X\oplus\kg_1$, we have
$[\kg,X]=[\kg_1,X]$. We choose functions $(\ga_1,\dots,\ga_r)$ on $G$ whose differentials
at the unit $e$ of $G$ are the equations of $\kg_1$ in $\kg$ and define $A=\ensemble{g\in
G}{\ga_1(g)=\dots=\ga_r(g)}$. Then there is a Zarisky open subset $U$ of $A\x \kg^Y$
containing $(e,0)$ which is smooth and such that the map $\gn:U\to \kg$ is \'etale.

Let $(s_1,\dots,s_r)$ be the coordinates of $\kg^Y$ associated with the basis $(Y_1,\dots,Y_r)$,
they define functions $(s_1,\dots,s_r)$ on $U$ and $\gn(s^{-1}(0))$ is equal to $S$. Let
$\gh_0=\gn^*E$ on $U$ ($E$ be the Euler vector field of $\kg$). A standard calculation \cite[Part
I, \S 5.6]{VARAD2}, shows that $\gh_0(s_i)=(\gl_i/2+1)s_i$ hence the map $F_0:\kg^Y\to \kg$ defined
by $F_0(Z)=X+Z$ is quasi-homogeneous if the weights on $\kg^Y$ are $m'_i=(\gl_i/2+1)$ for
$i=1,\dots,r$ and the weights on $\kg$ are $(1,\dots,1)$. The map $F:\kg^Y\to W$ defined by
$F(Z)=\gy(X+Z)$ is thus quasi-homogeneous with the weights $(m'_1,\dots,m'_r)$ on $\kg^Y$ and
$(m_1,\dots,m_l)$ on $W$.

Now, we apply corollary \ref{cor:invbf} to $X'=U$, $X=W$, $f=\gy_\circ\gn$, $W'=\kg^Y$,
$\gf'$ the projection $U\to\kg^Y$ and $F:\kg^Y\to W=k^l$ given by $F(Z)=\gy(X+Z)$. This
gives the $b$-function for $\gn^+\NO$ and thus for $\NO$ by lemma \ref{lem:etale}.
\end{proof}

Let us now consider the non-nilpotent strata of the stratification of $\kg$
(\S\ref{sec:stratification}):

\begin{prop}\label{prop:strati}
The module $\NO$ admits a tame quasi-$b$-function $b_S$ along each stratum $S$.

More precisely, if the stratum is $S_{(P,\kO)}$ according to definition \ref{def:strat} and $\kqp$
the associated semi-simple Lie subalgebra of $\kg$, then

a) $b_S$ depends only on $P$ and its roots are integers greater or equal to $-(m-k)/2$ where $m$ is
the dimension of $\kqp$ and $k$ its rank.

b) The total weight of $b_S$ is equal to $(m+r)/2$ where $r$ is the codimension of $\kO$ in $\kqp$.

In particular, on the stratum of codimension $1$ in $\kg$, the roots of the usual $b$-function of
$\NO$ are half integers greater or equal to $-1/2$.
\end{prop}

\begin{proof}
We fix a Cartan subalgebra of $\kg$ and a subset $P$ of roots with the notations of
\S\ref{sec:stratification}. This define a semi-simple algebra $\kqp$ to which are
associated the maps $\gy_P:\kqp\to W_P$ and $\gf_P:\khp\to W_P$. Here $W_P$ is a vector
space of dimension the rank of $\kqp$. The Cartan subalgebra $\kh$ splits into the direct
sum $\kh=\khp\oplus\khpo$ and this define a map $\gf'_P=\gf_P\otimes 1:\kh\to
W=k^l\oplus\khpo$.

Let $S\in\khpo$, we know from section \ref{sec:stratification} that
$\kg^S=\khpo\oplus\kqp$ and as $S$ is semisimple we have $\kg=[\kg,S]\oplus\kg^S$. The
map $\gn:G\times\kg^S\to\kg$ defined by $\gn(g,Z)=g.(Z+S)$ is thus a submersion. Let
$(\ga_1,\dots,\ga_r)$ be functions on $G$ whose differentials at $e$ are the equations of
$[\kg,S]$ in $\kg$ and define $A=\ensemble{g\in G}{\ga_1(g)=\dots=\ga_r(g)}$. Then there
is a Zarisky open subset $U$ of $A\x \kg^S$ containing $(e,0)$ which is smooth and such
that the map $\gn:U\to \kg$ is \'etale. Let $\gy':U\to W$ be defined as the composition
of the canonical projection  $A\x \kg^S\to\kg^S$ and of $\gy_P$.

Now we follow the proof of proposition \ref{prop:Euler} with the same notations. Applying the
second part of proposition \ref{prop:bf} to $L=\{0\}\x\khpp$, we find that $\Nl$ admits a
monodromic $b$-function along $L$ which is equal to
$b_0(\gvtp)=\gvtp(\gvtp-1)\dots(\gvtp-N'+1)(\gvtp-N')$ where $\gvtp$ is the Euler vector field of
$\khp$ and $N'$ is less or equal to $N=\frac{n-l}2$ with $n$ the dimension of $\kg$. This means
that there exists $l$ differential operators $R,A_1,\dots,A_l$ on $\kh$ that we may assume
invariant under $\kw_P$, with $R$ of order $-1$ for the $V$filtration associated with $L$ such that
$b_0(\gvtp)+R=A_1(z,D_z)q_1(D_z)+\dots+A_l(z,D_z)q_l(D_z)$. If $\gl=0$ we have $R=0$.

As these operators are invariant under $\kw_P$ hence under $\gf'$ we may apply  $\gf'_*$ and
$\gg^{-1}$ and find an equation
$b_0(\gg^{-1}(\gh))+\gg^{-1}\gf_*(R)=B_1\gy_*(Q_1)+\dots+B_l\gy_*(Q_l)$ with $B_1,\dots,B_l$ in
$\CaD(W)$ and $\gh=\gf_*(\gvt)$. In the coordinates of $W$ defined by the isomorphism
$\gf':k^l\oplus\khpp\to W$,  the vector $\gh$ is equal to $\sum n_it_iD_{t_i}$ where the $n_i$ are
the primitive degrees of $\khp$, it is associated with the manifold $L'=\gf'(\{0\}\oplus\khpp)$.

As $\gD$ is the product of $\gD_P$ by a function which does not vanish on a neighborhood of $L$,
the function $d$ which defines the morphism $\gg$ is the product of the corresponding function
$d_P$ associated with $\khp$ by a function $\gro$ which does not vanish in a neighborhood of $L'$.
So we have
\[\gg^{-1}(\gh)=\gro^{-1/2}d_P^{-1/2}\gh d_P^{1/2}\gro^{1/2}
=\gro^{-1/2}\left(\gh+N_P\right)\gro^{1/2}=\left(\gh+N_P\right)+a\] where $N_P$ is
$(m-k)/2$ ($m$ is the dimension of $\kqp$, $k$ its rank) and $a$ is a function which
vanishes on $L'$ hence of order at most $-1$ for the $V^\gh$-filtration. The operator
$\gg^{-1}\gf_*(R)$ is also of order $-1$ for the $V^\gh$-filtration.

We have proved that $\Nl$ admits a $b(\gh)$-function along $L'$ which is equal to
$b(T)=(T-N_P)\dots(T-N_P+N)$. The end of the proof is the same than to the proof of
proposition \ref{prop:Euler2}.
\end{proof}

Proposition \ref{prop:strati} shows that $\NO$ is tame. To prove theorem \ref{thm:second} we have
still to prove that it is conic. This come from the fact that the singular support of $\NO$ is
conic for the Euler vector field of $\kg$ and the vector fields associated with the strata are
equal to this Euler vector field modulo vector fields tangent to the orbits.

\subsection{Isomorphism with the inverse image.}

We recall that $\NO$ is the submodule of $\gy^+\Nl$ generated by $u_{\kg\to W}$, it is
the image of the morphism $\MF\to\gy^+\Nl$.

\begin{theo}\label{thm:iso}
The canonical morphism $\MF\to\NO$ is an isomorphism.
\end{theo}

\begin{proof}
\textbf{1st step:}From semi-simple Lie algebras to reductive algebras.

Assume that the result has been proved for semi-simple Lie algebras and let $\kg$ be a reductive
algebra, direct sum of its center and a semisimple Lie algebra. By induction, we may assume that
$\kg=\kc\oplus\kg'$ with $\kc$ subspace of the center of dimension $1$ and $\kg'$ reductive Lie
algebra for which the result has been proved.

Let $t$ a coordinate of $\kc$ and $\gt$ the corresponding coordinate of the dual space
$\kc^*$. By the hypothesis, there is a differential operator in $F$ whose principal
symbol is equal to some power $\gt^q$. This means that $\kg'=\{t=0\}$ is non
characteristic for $\MF$. Let $\CK$ be the kernel of $\MF\to\NO$. We have an exact
sequence $\exacte{\CK}{\MF}{\NO}$ of non characteristic $\Dg$-modules. As the inverse
image is an exact functor in the non characteristic case, this gives an exact sequence
$\exacte{\CK/t\CK}{\MF/t\MF}{\NO/t\NO}$. If we prove that $\CK/t\CK=0$, we will have
$\CK=0$ (as $\CK$ is non characteristic).

So, we have to prove that $\MF/t\MF\to\NO/t\NO$ is injective. Here we use the same proof than in
\cite[Lemma 2.2.3.]{LBY}. In fact, as $\CaD_{\kg'}$-module $\MF/t\MF$ is generated by the classes
of $1,D_t,\dots,D_t^{q-1}$ and the submodule generated by $D_t^{q-1}$ is a module on $\kg'$ of the
same type than $\MF$ for which the theorem is true. Then we consider the quotient of $\MF/t\MF$ by
the module generated by $D_t^{q-1}$ and argue by induction.

\textbf{2nd step:}The result is true at points $X\in\kg$ whose semi-simple part is non zero.

By the first step, we may assume that $\kg$ is semisimple. Let $S$ be a non zero semisimple element
of $\kg$, $\kg^S$ its centralizer and $G^S$ the corresponding group. The spaces $\CO(\kg)^G$ and
$\CO(\kg^S)^{G^S}$ are isomorphic hence the space $W_S$ associated with $\kg^S$ is equal to $W$ and
the map $\gy_S:\kg^S\to W$ is the restriction of $\gy:\kg\to W$. Thus the sheaf of differential
operators on $\kg^S$ invariant under the action of $G^S$ is isomorphic to $\CaD_\kg^G$.

By induction on the dimension of $\kg$, we may assume that the theorem is true for $\kg^S$ hence
that the morphism $\CM_F^S\to\gy_S^+\Nl$ is injective. Here $\CM_F^S$ is the $\CaD_{\kg^S}$ module
associated with $F$ and $\CN_F$ the quotient of $\DV$ by the ideal generated by $F$. By definition,
the germ at $S$ of $\gy^+\CN_F$ is
$(\gy^+\CN_F)_S=\CO_{\kg,S}\otimes_{\CO_{\kg^S,S}}(\gy_S^+\Nl)_S$. On the other hand, we have
$(\Dg/\Dg\gt(\kg))_S=\CO_{\kg,S}\otimes_{\CO_{\kg^S,S}}(\CaD_{\kg^S}/\CaD_{\kg^S}\gt(\kg^S))_S$
hence $\CM_{F,S}=\CO_{\kg,S}\otimes_{\CO_{\kg^S,S}}\CM_{F,S}^S$. The morphism $\MF\to\gy^+\CN_F$ is
thus injective at the point $S$ hence at all the orbits whose closure contains $S$ that is in
particular at all points $X$ whose semisimple part in the Jordan decomposition is $S$.

\textbf{3nd step:}The case of nilpotent orbits.

Let $\CK$ be the kernel of $\MF\to\NO$. By the second step, we may assume that the theorem is true
at all non nilpotent points of $\kg$ that is that $\CK$ is supported by the nilpotent cone. Let
$\CK(\kg)^G$ be the set of global sections of $\CK$ invariant under $G$, we get an exact sequence
$$\exacte{\CK(\kg)^G}{\CM(\kg)^G}{\CM_0(\kg)^G}$$
and by \cite[lemma 3.2.]{HUNZ} we have $\CM(\kg)^G=\CM_0(\kg)^G=\Nl$ hence $\CK(\kg)^G=0$. Then
$\CK=0$ by \cite[lemma 3.2.]{LEVASS2}.
\end{proof}

Remark that the third step is also a consequence of the property (5) of the
Harish-Chandra morphism which has been proved by Levasseur-Stafford \cite{LEVASS3}.

\providecommand{\bysame}{\leavevmode ---\ } \providecommand{\og}{``} \providecommand{\fg}{''}
\providecommand{\smfandname}{\&} \providecommand{\smfedsname}{\'eds.}
\providecommand{\smfedname}{\'ed.} \providecommand{\smfmastersthesisname}{M\'emoire}
\providecommand{\smfphdthesisname}{Th\`ese}

\enddocument

\end